\documentclass[reqno,10pt]{amsart}
\usepackage{amsthm}
\usepackage{amsmath}
\usepackage{amssymb}
\usepackage{amsfonts}
\usepackage{mathrsfs}

\usepackage[utf8]{inputenc}
\usepackage{latexsym}
\usepackage{verbatim}

\usepackage[all]{xy}
\usepackage[bookmarksnumbered,colorlinks,urlcolor=blue]{hyperref}
\usepackage{graphicx}
\usepackage{graphics}
\usepackage{float}
\usepackage{enumerate}
\def\bb#1\eb{\textcolor{blue}
{#1}} %
\def\br#1\er{\textcolor{red}
{#1}} %
   \def\br#1\er{\textcolor{red}{#1}} %
      \def\bb#1\eb{\textcolor{blue}{#1}} %

\title[]{Geodesic connectedness of affine manifolds}

\author[I.P. Costa e Silva]{Ivan P. Costa e Silva}
\address{Department of Mathematics, Universidade Federal de Santa Catarina,
	\hfill\break\indent 88.040-900
	Florianópolis-SC, Brazil.}
\email{pontual.ivan@ufsc.br}

\author[J.L. Flores]{José L. Flores}
\address{Departamento de \'Algebra, Geometr\'{\i}a y Topolog\'{\i}a,  Universidad de M\'alaga
\hfill\break\indent
Facultad de Ciencias, Campus Universitario de Teatinos,
\hfill\break\indent 29071 M\'alaga, Spain}
\email{floresj@uma.es}


\begin{document}
\newtheorem{thm}{Theorem}[section]
\newtheorem{prop}[thm]{Proposition}
\newtheorem{lemma}[thm]{Lemma}
\newtheorem{cor}[thm]{Corollary}
\theoremstyle{definition}
\newtheorem{defi}[thm]{Definition}
\newtheorem{notation}[thm]{Notation}
\newtheorem{exe}[thm]{Example}
\newtheorem{conj}[thm]{Conjecture}
\newtheorem{prob}[thm]{Problem}
\newtheorem{rem}[thm]{Remark}

\newtheorem{conv}[thm]{Convention}
\newtheorem{crit}[thm]{Criterion}
\newtheorem{claim}[thm]{Claim}

\newcommand{\ben}{\begin{enumerate}}
\newcommand{\een}{\end{enumerate}}
\newcommand{\qcd}{\begin{flushright} $\Box$ \end{flushright}}
\newcommand{\bit}{\begin{itemize}}
\newcommand{\eit}{\end{itemize}}

\begin{abstract}
We discuss new sufficient conditions under which an affine manifold $(M,\nabla)$ is geodesically connected. These conditions are shown to be essentially weaker than those discussed in groundbreaking work by Beem and Parker and in recent work by Alexander and Karr, with the added advantage that they yield an elementary proof of the main result. 
\end{abstract}

\maketitle

\vspace*{-.5cm}

\section{Introduction}\label{section1}

The geodesic connectedness of a (connected) Riemannian manifold $(M,g)$ is of course directly linked with geodesic completeness via the Hopf-Rinow theorem, in which the existence of minimal geodesics connecting any two points of $M$ is established. Even in the absence of geodesic completeness, the geodesic connectedness of Riemannian manifolds is fairly well-understood \cite{miguelreview,rosselathesis} 
By contrast, it has long been known that for indefinite semi-Riemannian manifolds, where no analogue of the Hopf-Rinow theorem exists, the situation is much subtler. 

A famous example by Bates \cite{BATES} has shown that even {\em complete and compact} affine manifolds may fail to be geodesically connected. Even if one only considers the more restricted (but important) class of {\em Lorentzian} manifolds, it is well-known that de Sitter and anti-de Sitter spaces are examples of geodesically complete, maximally symmetric Lorentz manifolds which are not geodesically connected. The underlying manifolds in these examples, however, are not compact. Yet, compact Lorentz tori can still be found which also fail to be geodesically connected. These examples, in turn, are not geodesically complete\footnote{It is well-known, again due to the absence of a Hopf-Rinow-like result, that compactness and geodesic completeness are unrelated in Lorentzian manifolds (cf. Example 7.17 in \cite{oneill}).}. Indeed, to the best of our knowledge it remains an open problem to ascertain whether a connected Lorentzian manifold which is {\em both} compact {\em and} geodesically complete is geodesically connected. 

On the flip side, powerful variational tools have been developed starting with the seminal work of Benci and Fortunato \cite{BF1,BF2} to tackle the problem of geodesic connectedness in important special cases of Lorentzian manifolds. These special cases not only include stationary (i.e., endowed with a complete timelike Killing vector field), but also time-dependent orthogonal splitting Lorentzian manifolds. But an important limitation of these methods ultimately arises from the fact that the natural energy functional on the space of curves does not, in general, satisfy any natural compactness property - such as the Palais-Smale condition - for indefinite metrics. Thus, only certain quite symmetric classes of semi-Riemannian manifolds can be tackled through variational methods (see, for example, the review works \cite{CaSa} and references therein for a vast number of examples that can). 

Other techniques which have proven to be fruitful for the study of geodesic connectedness are: methods from group theory \cite{CM}, topological techniques based on Brouwer's degree (see \cite{F} and references therein), and more direct but restricted methods based on a (partial) integration of the geodesic equations (see, for instance, \cite{CS}).

In an entirely different vein, there are the works of Beem \cite{Beem1} and Beem $\&$ Parker \cite{Beem_Parker1,Beem_Parker2}, which do not use the previous techniques at all, but rather give two quite general geometric conditions, {\em disprisonment} and {\em pseudoconvexity}, which, unlike variational methods, apply to families of geodesics in the larger class of {\em affine manifolds}. Specifically they prove \cite[Prop.\ 5]{Beem_Parker2}:
\begin{thm}\label{beemParkerthm}
Let $M=(M,\nabla)$ be a connected disprisoning and pseudoconvex affine manifold without conjugate points. Then  $M$ is geodesically connected.
\end{thm}
For the purposes of this note, we may divide the geometric hypotheses in this theorem in two groups:
\begin{itemize}
\item[1)] disprisonment and pseudoconvexity, 
\item[2)] absence of conjugate points.
\end{itemize}
We shall recall the precise definitions of hypotheses $(1)$ in the next section. Suffice it to say here that what these assumptions really amount to is ensuring for each $p \in M$ that the exponential map $\exp_p:\mathcal{D}_p \subset T_pM \rightarrow M$ is a {\em proper} map (cf. \cite[Lemma \ 4]{Beem_Parker2}), i.e., the preimage of any compact set in $\mathcal{D}_p$ is compact. In particular, it has closed image. Then, hypothesis $(2)$ implies it is also a local diffeomorphism, and hence an open map, whence the conclusion follows by the connectedness of $M$. 

In this paper we separately weaken {\em both} $(1)$ and $(2)$. Instead of $(1)$, we introduce (what we show to be) a strictly weaker condition on exponential maps, called on that account {\em weak properness}, which, unlike properness, applies to {\em compact} manifolds as well, and moreover holds in any complete Riemannian manifold, whether compact or not. In lieu of $(2)$, we introduce the requirement that for any $p \in M$, $(i)$ the set ${\rm Conj}(p)$ of conjugate points to $p$ is closed and $(ii)$ $M\setminus {\rm Conj}(p)$ is connected. 
We call connected affine manifolds satisfying the latter condition
{\em weakly Wierdesehen manifolds}, or {\em WW manifolds} for short. Clearly 
these conditions are trivial if the manifold has no conjugate points, but allow for many more examples otherwise\footnote{We remind the reader that a connected complete Riemannian manifold is said to be a {\em Wierdersehen} manifold if ${\rm Conj}(p)$ is just one point for every $p \in M$. The importance of these manifolds arises due to the so-called (original) {\em Blaschke conjecture}, which posited that any Wierdersehen manifolds of $dim\ M=2$ is isometric to a round sphere $\mathbb{S}^2$. This has been proven by Green \cite{green}, and the analogue statement in higher dimensions has been proven by Berger, Kazdan, Weinstein and Yang \cite{Ber,Kaz,Wein,Y}.}, such as (if $dim\ M >1$) when each $p \in M$ always has at most a {\em discrete} set of conjugate points. 

Our main theorem is the following. 
\begin{thm}\label{mainthm}
Let $M=(M,\nabla)$ be a WW affine manifold. If $\exp_p:\mathcal{D}_p \subset T_pM \rightarrow M$ is weakly proper for each $p \in M$, then $M$ is geodesically connected. 
\end{thm}

\begin{rem} It is not clear to the authors to what extent the assumption that the manifold be WW is actually necessary, or whether it is just an artifact of the limitations of the approach adopted here. In any case, the following statements will be easily verifiable after one gets acquainted with the proof of Theorem \ref{mainthm}:
\begin{itemize}	
\item[(1)] If the conditions (i) and (ii) of the notion of WW manifold are replaced at some point $p$ by $M\setminus \overline{{\rm Conj}(p)}$ being connected, then $p$ can be geodesically connected with ``almost any'' point of $M$ (eventually unreachable-from-$p$ points must be contained in the zero-meaure set $\overline{{\rm Conj}(p)}\setminus {\rm Conj}(p)$).  
\item[(2)] If the conditions (i) or (ii) are not satisfied at some point $p$, then $p$ can be geodesically connected with any point of any connected component of $M\setminus \overline{{\rm Conj}(p)}$ whose intersection with some normal neighborhood $\mathcal{V}$ of $p$ is non-empty.
\end{itemize}
\end{rem}


More recently, Alexander and Karr \cite{AK,K} have also given a geodesic connectedness result on semi-Riemannian manifolds possessing a suitable {\em convex function}: 
\begin{thm}\label{AlexKarr}
Let $M$ be a [connected] null-disprisoning semi-Riemannian manifold. Suppose
$M$ supports a proper, nonnegative convex function $f : M \rightarrow \mathbb{R}$ whose critical set is a minimum point. If there is no non-constant complete geodesic on which $f$ is constant (for example, if $f$ is strictly convex), then $M$ is geodesically connected.
\end{thm}
We also show here that {\em the existence of such a convex function $f$ also implies that exponential maps are all proper}. In fact, we show that under the assumptions in Theorem \ref{AlexKarr} conditions $(1)$ hold, so in a way the existence of $f$ is actually a {\em stronger} condition than Beem and Parker's. On the other hand, no assumption on conjugate points is made in Theorem \ref{AlexKarr}, so no version of condition $(2)$ needs to be used therein. But while our main result restricts to WW manifolds, the proof we present here, after suitable definitions are given, will be quite elementary, apart from being applicable to more general affine manifolds than just semi-Riemannian ones.

It should be noted that the central argument provided in this paper is very general, and can be cast in a wider context. In order to stress this fact, we deduce Theorem \ref{mainthm} via the proof of the more general Theorem \ref{mmain}.

The rest of this note is divided in two sections: in section \ref{s1} we introduce the necessary definitions and ancillary results we shall need, and the proof of Theorem \ref{mainthm} is given in section \ref{s2}.

\section{Technical preliminaries}\label{s1}

We begin by recalling (cf. \cite{BeemBook,Beem_Parker1,Beem_Parker2}) that a collection $\mathcal{C}$ of geodesics on an affine manifold\footnote{Here and hereafter, {\em smooth} means $C^{\infty}$ and {\em manifolds} for us always means {\em real, smooth, finite-dimensional, Hausdorff, second-countable} manifolds.} $(M,\nabla)$ is 
\begin{itemize}
\item[a)] {\em pseudoconvex} if for any compact set $K\subset M$ there exists a compact set $K^{\ast}\subset M$ such that any geodesic in $\mathcal{C}$  with endpoints in $K$ is entirely contained in $K^{\ast}$;
\item[b)] {\em imprisoned} if some maximal geodesic in $\mathcal{C}$ is entirely contained in some compact subset of $M$. Otherwise $\mathcal{C}$ is said to be {\em disprisoned}. 
\end{itemize}
If we take $\mathcal{C}$ to be the collections of {\em all} geodesics in $(M,\nabla)$ and that class is pseudoconvex and disprisoned, then $(M,\nabla)$ itself is said to be a ({\em geodesically}) {\em pseudoconvex} and {\em disprisoning} affine manifold. Another case of interest below will be when $\nabla$ is the Levi-Civita connection of a semi-Riemannian metric $g$ and $\mathcal{C}$ is the collection of {\em null} (or {\em lightlike}) {\em geodesics}. If that collection is pseudoconvex/disprisoned, then $(M,g)$ itself is said to be {\em null-pseudoconvex}/{\em null-disprisoning}, respectively.

\bigskip
Recall also that a map $F:X\rightarrow Y$ between topological spaces $X,Y$ is {\em proper} if $F^{-1}(K)$ is compact in $X$ whenever $K\subset Y$ is compact. Of course, if $Y$ is compact but $X$ is not, then $F$ can never be proper. Given a semi-Riemannian manifold $(M,g)$, a smooth function $f:M\rightarrow \mathbb{R}$ is {\em convex} [resp. {\em strictly convex}] if 
its Hessian ${\rm Hess}_g f$ is positive-semidefinite [resp. positive-definite]. 

\begin{prop}\label{prop1}
Let $(M,\nabla)$ be an affine manifold.  
\begin{itemize}
\item[i)] If  $(M,\nabla)$ is both pseudoconvex and disprisoning, then $\forall p \in M$ the exponential map $\exp_p:\mathcal{D}_p\subset T_pM\rightarrow M$ is a proper map.
\item[ii)] Suppose $\nabla$ is the Levi-Civita connection of a semi-Riemannian metric $g$ on $M$, and that $(M,g)$ is null-disprisoning and possesses a smooth proper convex function $f:M\rightarrow \mathbb{R}$ satisfying in addition all the conditions in Thm. \ref{AlexKarr}. Then $(M,g)$ is both pseudoconvex and disprisoning.
\end{itemize}
\end{prop}
\noindent {\em Proof.} $(i)$ As indicated in the Introduction, this follows rather easily from the proof of \cite[Lemma \ 4]{Beem_Parker2}.

$(ii)$ That $(M,g)$ is disprisoning follows directly from \cite[Lemma \ 3.3]{AK}, so we only need to show it is pseudoconvex. Let $\gamma:I\subset \mathbb{R}\rightarrow M$ be any geodesic (here $I$ denotes any nonempty interval). Then:
\begin{eqnarray}
{\rm Hess}_g f(\gamma ', \gamma ') &=& g(\nabla _{\gamma '}\nabla f \circ \gamma, \gamma ') \nonumber \\
&=& (g(\nabla f\circ \gamma ,\gamma '))'- g(\nabla f \circ \gamma , \nabla _{\gamma '}\gamma ') \nonumber \\
&\equiv & (f\circ \gamma)'', \nonumber
\end{eqnarray}
whence we conclude, by the convexity of $f$, that 
\begin{equation}\label{secondder1}
(f\circ \gamma)''\geq 0.
\end{equation}
Given $a,b \in I$ with $a<b$, elementary Calculus then yields by (\ref{secondder1}), since we are assuming $f\geq 0$, that 
\begin{equation}\label{secondder2}
0\leq f\circ \gamma |_{[a,b]}\leq \mbox{max}\{f(\gamma(a)),f(\gamma(b))\}.
\end{equation}
Let $(K_n)_{n \in \mathbb{N}}$ be a sequence of compact subsets of $M$ such that 
\[
K_{n+1}\subset \mbox{int}\,K_n, \quad \bigcup_{n \in \mathbb{N}}K_n = M.
\]
Assume now that $(M,g)$ is not pseudoconvex. Then we may pick, for some compact $K_0\subset M$, a sequence of geodesics $\gamma _n:[a_n,b_n]\rightarrow M$ such that for each $n$, $\gamma _n(a_n),\gamma _n(b_n)\in K_0$, but there exists $t_n\in (a_n,b_n)$ such that $\gamma_n(t_n)\notin K_n$. 
Let $A>0$ be the maximum value $f$ takes on $K_0$. Thus, (\ref{secondder2}) gives
\[
0\leq f\circ \gamma_n \leq A, \quad \forall n \in \mathbb{N},
\]
i.e., $\gamma _n[a_n,b_n]\subset f^{-1}([0,A]), \forall n \in \mathbb{N}$. But since $f$ is proper, $f^{-1}([0,A])$ is compact, so for large enough $n_0$,
\[
\gamma_{n_0}[a_{n_0},b_{n_0}] \subset f^{-1}([0,A]) \subset \mbox{int}\,K_{n_0}, 
\]
in contradiction with $\gamma(t_{n_0})\notin K_{n_0}$. 
\qcd

Given a topological space $X$, we say that a continuous map $\alpha:[a,b)\rightarrow X$ with $-\infty< a<b\leq +\infty $ in the extended real numbers $\overline{\mathbb{R}}$ is {\em right-extendible} if it admits a continuous extension $\tilde{\alpha}:[a,b]\rightarrow X$, or equivalently, if it possesses a {\em right-endpoint}, i.e. some $x\in X$ such that for any open set $U\ni x$ there exists $a<c<b$ such that $\alpha[c,b)\subset U$. (If $X$ is Hausdorff, then such a right-endpoint or continuous extension is unique if it exists.) If $\alpha$ is not right-extendible, then it is said to be {\em right-inextendible}.

A key definition for us will be the following. 
\begin{defi}\label{wpdef}
Let $M_1,M_2$ be smooth manifolds, and let $F:M_1\rightarrow M_2$ be a smooth map. We say that $F$ is {\em weakly proper at $p\in M_1$} if any piecewise smooth regular curve $\alpha:[a,b)\rightarrow M_1$, $-\infty< a<b\leq +\infty $, with $\alpha(a)=p$ for which $F\circ \alpha$ is right-extendible in $M_2$ has its image contained in some compact set in $M_1$. $F$ is {\em weakly proper} if it is weakly proper at any $p \in M$. 
\end{defi}

\begin{rem}\label{rem1}
Clearly, if $F:M_1\rightarrow M_2$ is weakly proper at some $p \in M_1$ and $M_1$ is connected, then $F$ is weakly proper.
\end{rem}

We first have:
\begin{prop}\label{prop2}
Let $M_1,M_2$ be smooth manifolds. A smooth proper map $F:M_1\rightarrow M_2$ is weakly proper.
\end{prop}
\noindent {\em Proof.} Let $p \in M_1$. Given any piecewise smooth curve $\alpha:[a,b)\rightarrow M_1$ with $\alpha(a)=p$, if $F\circ \alpha$ is right-extendible in $M_2$, then in particular $F\circ \alpha[a,b)$ is precompact. Therefore
\[
\alpha[a,b) \subset F^{-1}( \overline{F\circ \alpha[a,b)}),
\]
and the set on the right is compact if $F$ is proper. 
\qcd
Props. \ref{prop1} and \ref{prop2} now immediately yield:
\begin{cor}\label{cor1}
If  $(M,\nabla)$ is both pseudoconvex and disprisoning, then $\forall p \in M$ the exponential map $\exp_p:\mathcal{D}_p\subset T_pM\rightarrow M$ is weakly proper. 
\end{cor}
\qcd 
If  $(M,\nabla)$, however, is both pseudoconvex and disprisoning, then we have seen that $M$ can never be compact, while there is no such restriction when only weak properness holds for its exponential maps. 
\begin{prop}\label{prop3}
If $(M,g)$ is a complete Riemannian manifold, then $\forall p \in M$ the exponential map $\exp_p:T_pM\rightarrow M$ is weakly proper. 
\end{prop}
\noindent {\em Proof.} Fix $p \in M$. By Remark \ref{rem1}, it is enough to show that $\exp_p$ is weakly proper at the origin $0_p\in T_pM$. Let $\overline{\alpha}:[a,b)\rightarrow T_pM$ be a piecewise smooth curve with $\overline{\alpha}(a)=0_p$. Let $J$ be the set of $t\in [a,b)$ with $\overline{\alpha}(t)\neq 0_p$, and for any $t\in J$, put
\[
r(t):= |\overline{\alpha}(t)|_{g_p} \mbox{  and  }  w(t):= \frac{\overline{\alpha}(t)}{r(t)}.
\]
Let 
\[
\varphi:(r,s) \in (0,+\infty) \times J \mapsto \exp_p(r\cdot w(s)) \in M. 
\]
Put $\alpha := \exp_p\circ \overline{\alpha}$, so that we can write 
\[
\alpha '(t) = r'(t)\frac{\partial \varphi}{\partial r}(r(t),t) + \frac{\partial \varphi}{\partial t}(r(t),t), \quad \forall t \in J.
\]
By Gauss' Lemma we conclude that $\langle \frac{\partial \varphi}{\partial r},\frac{\partial \varphi}{\partial r} \rangle =1$ and $\langle \frac{\partial \varphi}{\partial r},\frac{\partial \varphi}{\partial t} \rangle =0$, so that
\begin{equation}\label{refhere}
g(\alpha '(t),\alpha '(t))= (r'(t))^2 + \left|\frac{\partial \varphi}{\partial t}(r(t),t)\right| ^2 \geq (r'(t))^2, \quad \forall t \in J.
\end{equation}
Assume $\alpha:[a,b)\rightarrow M$ is right-extendible. Since it is piecewise smooth and $g$ is complete we must have 
\[
\ell:= L_g(\alpha) = \int _a^b |\alpha '(s)|_g ds  <+\infty.
\]
Assume, by way of contradiction, that $\overline{\alpha}[a,b)$ is not precompact. Then, it is also not bounded, so there exists an increasing sequence $(t_n)\subset J$ with  such that 
\[
r(t_n)> n, \quad \forall n\in \mathbb{N}.
\]
Hence, by integrating Eq. (\ref{refhere}), we get 
\[
\ell = \int _a^b |\alpha '(s)|_g ds \geq n, \quad \forall n\in \mathbb{N},
\]
an absurd. We conclude that $\overline{\alpha}[a,b)$ is indeed precompact, and hence $\exp_p$ is weakly proper. 
\qcd

We also provide here a sufficient condition for weak properness in the Lorentzian case. 
\begin{prop}\label{prop4}
	If $(M,g)$ is a complete Lorentzian manifold with a parallel timelike vector field\footnote{As proven in \cite{GM}, this condition is equivalent to precompactness of the holonomy group of $(M,g)$.}, then $\forall p \in M$ the exponential map $\exp_p:T_pM\rightarrow M$ is weakly proper. 
\end{prop}
\noindent {\em Proof.} Let $V\in \mathfrak{X}(M)$ be a timelike parallel field. Since $V$ is parallel, $g(V,V)$ is constant, and thus we assume without loss of generality that $g(V,V)=-1$. 
Fix a {\em Riemannian} metric $h$ on $M$ as follows. For any $X,Y \in \mathfrak{X}(M)$, define
\[
h(X,Y):= 2g(V,X)g(V,Y) + g(X,Y).
\]
Note that $h(V,V) =1$ and that for any $X \in \mathfrak{X}(M)$,
\begin{equation}\label{eqclear}
h(V,X)=2g(X,V)g(V,V) + g(X,V) = -2g(V,X) + g(V,X) = - g(V,X).
\end{equation}
(In particular, $X$ is $h$-orthogonal to $V$ if and only if it is $g$-orthogonal to $V$.) As a convenient notation, we adopt $\langle \, . \, , \, . \, \rangle := g$. Following the standard notation in (semi-)Riemannian geometry, for any $v\in TM$ we put
\[
|v|^2:=  \langle v,v\rangle,
\]
which can have any sign, but the {\em norm} $|v| :=\sqrt{|\langle v,v\rangle|}$ is always nonnegative. The corresponding quantities for $h$ are of course always nonnegative, and any $h$-related quantities will be explicitly accompanied by a subscript ``$h$'' to avoid confusion.  (So, for example, $|v|_h$ means the $h$-norm of $v\in TM$.)

Assume by contradiction that $\exp_p:T_pM\rightarrow M$ is not weakly proper for some $p\in M$. Then, there exists a piecewise smooth curve $\overline{\gamma}:[0,\ell)\rightarrow T_p M$ whose image eventually leaves any compact subset in $T_pM$, but such that $\gamma=\exp_p\circ \overline{\gamma}$ is right-extendible. Due to Remark \ref{rem1}, we can also assume that $\overline{\gamma}(0)=0$. 

Right-extendibility means that the image $\gamma[0,\ell)$ is precompact in $M$, and that its $h$-length is finite; in other words, $L_h(\gamma)<\infty$. We can assume, without loss of generality, that $\gamma$ is $h$-arc length parametrized (and thus $|\gamma'|_h =1$ and $\ell(=L_h(\gamma))<\infty$). We do so for the remaining of the proof. 

Now, consider the $g$-geodesic variation given by 
\[
\sigma:(s,t) \in [0,1]\times [0,\ell) \rightarrow \exp_p(s\cdot \overline{\gamma}(t)) \in M .
\]
(See, e.g., p. 122 of \cite{oneill} for background on manipulating with variations.) Then 
\begin{eqnarray}
\frac{\partial \sigma }{\partial s}(s,t) &=& ( d \exp _p)_{s\cdot \overline{\gamma}(t)}(\overline{\gamma}(t)_{s\cdot \overline{\gamma}(t)}), \label{eq16.1} \\
\frac{\partial \sigma }{\partial t}(s,t) &=& ( d\exp _p)_{s\cdot \overline{\gamma}(t)}(s\cdot \overline{\gamma}'(t)_{s\cdot \overline{\gamma}(t)}), \label{eq16.2}\\
\frac{\partial ^2 \sigma }{\partial s \partial t}& \equiv &\frac{\partial ^2\sigma }{\partial t \partial s}, \label{eq16.3} \\
\frac{\partial (V\circ \sigma) }{\partial s}&=&\frac{\partial (V\circ \sigma) }{\partial t} \equiv 0, \label{eq16.4}
\end{eqnarray}
where equations (\ref{eq16.4}) are due to $V$ being parallel. 
It is a standard fact in Lorentzian geometry that since $V_p$ is timelike 
\[
T_pM= \mathbb{R}\cdot V_p \oplus (\mathbb{R}\cdot V_p)^{\perp},
\]
where the ``perp'' is {\em spacelike}, i.e., $g$ restricted to $(\mathbb{R}\cdot V_p)^{\perp}$ is positive definite.  Thus, for each $t \in [0,\ell)$ we can write 
\begin{equation}\label{eq17}
\overline{\gamma}(t)=a(t) \cdot V_p+ Z(t),
\end{equation}
with $Z$ a smooth vector field along $\overline{\gamma}$ which is $g$-normal (and hence $h$-normal) to $V_p$, i.e. $\langle V_p,Z \rangle_p = h_p(V_p,Z)\equiv 0$, so that $|Z|^2\geq 0$. Therefore, using Eqs. (\ref{eq16.1})-(\ref{eq16.4}), we compute
\begin{eqnarray}
I(t) & := & \int _0^1 ds \int _0^{t} d\lambda  \langle (V\circ \sigma)(s,\lambda), \frac{\partial ^2\sigma }{\partial t \partial s}(s,\lambda)\rangle \nonumber \\
&=& \int _0^1 ds \int _0^{t} d\lambda \frac{\partial \langle (V\circ \sigma), \frac{\partial \sigma }{\partial s}\rangle}{\partial t} (s,\lambda) \nonumber \\
&=& \int _0^1 ds [\langle (V\circ \sigma)(s,t), \frac{\partial \sigma }{ \partial s}(s,t)\rangle - \langle (V\circ \sigma)(s,0), \frac{\partial \sigma }{\partial s}(s,0)\rangle] \nonumber  \\
&\equiv & [\langle (V\circ \sigma)(0,t), \frac{\partial \sigma }{ \partial s}(0,t)\rangle - \langle (V\circ \sigma)(0,0), \frac{\partial \sigma }{\partial s}(0,0)\rangle]. 
\end{eqnarray}
(The last equality holds because $s\in [0,1] \mapsto \sigma (s,t)$ are geodesics, and thus their $s$-derivative and $V\circ \sigma(\, . \,, t)$ are parallel-transported along them.) A glance at Eqs. (\ref{eq16.1}) and (\ref{eq16.2}) show that 
\[
\frac{\partial \sigma }{\partial s}(0,0) \equiv 0 \mbox{  and  }  \frac{\partial \sigma }{\partial s}(t,0) \equiv \overline{\gamma}(t),\quad \forall t \in [0,\ell),
\]
whence Eq. (\ref{eq17}) yields 
\begin{equation}\label{eq18}
I(t) = \langle V_p,\overline{\gamma}(t)\rangle =-a(t).
\end{equation}
On the other hand, a similar manipulation with Eqs. (\ref{eq16.1})-(\ref{eq16.4}) also yield
\begin{eqnarray}
I(t)& = & \int _0^{t} d\lambda \int _0^1 ds\langle (V\circ \sigma)(s,\lambda), \frac{\partial ^2\sigma }{\partial t \partial s}(s,\lambda )\rangle \nonumber \\
&=&  \int _0^{t} d\lambda \int _0^1 ds \frac{\partial \langle (V\circ \sigma), \frac{\partial \sigma }{\partial t}\rangle}{\partial s} (s,\lambda ) \nonumber \\
&=&  \int _0^{t} d\lambda  \langle V\circ \gamma (\lambda),\gamma '(\lambda)\rangle. \label{addme}
\end{eqnarray}
Since $\gamma $ is $h$-arc-length-parametrized, and $V$ is $h$-unitary, we can use Eq. (\ref{eqclear}) and the Cauchy-Schwarz inequality to obtain
\begin{equation}\label{clearmore}
 |\langle V\circ \gamma (t),\gamma '(t)\rangle| = |h(V\circ \gamma (t),\gamma '(t))| \leq 1, \quad \forall t\in [0,\ell).
 \end{equation}
Now, substituting (\ref{clearmore}) in the last integral of (\ref{addme}) together with (\ref{eq18}) we get
\begin{equation}\label{clearevenmore}
 |a(t)|\leq \ell, \quad \forall t\in [0,\ell).
\end{equation}

We conclude that the function $ t\in [0,\ell) \mapsto a(t) \in \mathbb{R}$ is bounded. However,
\[
|\overline{\gamma}(t)|^2_h = |a(t)|^2+ |Z(t)|_h^2,
\]
and the boundedness of $a$ together with the assumed unboudedness of $ |\overline{\gamma}|_h$ allow us to conclude that 
\begin{equation}\label{eq19}
|Z|_h \equiv |Z|
\end{equation}
is not bounded in $T_pM$. Moreover, 
\[
|\overline{\gamma}|^2 = -|a|^2+ |Z|^2,
\]
therefore the set
\[
J:= \{t\in [0,\ell) \, : \, |\overline{\gamma}(t)|^2 >0 \} 
\]
is an open set in $[0,\ell)$ over which $ t\in J \mapsto |\overline{\gamma}(t)| \in \mathbb{R}$ is unbounded. 

Consider the smooth functions $r:J\rightarrow (0,+\infty)$ and $\overline{N}:J\rightarrow T_pM$ given by 
\[
r(t):=  |\overline{\gamma}(t)|, \mbox{   and   } \overline{N}(t):=  \overline{\gamma}(t) /|\overline{\gamma}(t)|, \quad \forall t \in J.
\]
Thus, we have 
\begin{eqnarray}
\overline{\gamma} &=& r \cdot \overline{N}, \label{eq20.1}\\
\overline{\gamma}' &=& r'\cdot \overline{N} + r\cdot \overline{N}' \Longrightarrow \gamma'= r'\cdot N +r\cdot Y, \label{eq20.2}
\end{eqnarray}
where we have defined smooth fields $N,Y\in \mathfrak{X}(\gamma |_J)$ by
\[
N(t):= (d\exp _p)_{\overline{\gamma}(t)}( \overline{N}(t)), \mbox{  and  } Y(t):= (d\exp _p)_{\overline{\gamma}(t)}( \overline{N}'(t)), \quad \forall t\in J.
\]
By Gauss' Lemma, $\langle N,N \rangle =1$ and $\langle N, Y\rangle =0$.  Thus, on the one hand we have, again using (\ref{eqclear}),
\begin{equation}\label{eq21}
|N|_h ^2 = h(N,N) = 2(h(V\circ \gamma, N))^2 + 1 \geq 1.
\end{equation}
On the other hand, $\forall t \in J$ we compute 
\begin{eqnarray}
\langle N(t), V_{\gamma(t)} \rangle &=&   \frac{1}{r(t)}\langle (d\exp _p)_{\overline{\gamma}(t)}( \overline{\gamma}(t)), V_{\gamma(t)} \rangle  \nonumber \\
&=& \frac{1}{r(t)}\langle \sigma _s(1,t),V _{\sigma(1,t)}\rangle  \nonumber \\
&=& \frac{1}{r(t)}\langle \sigma _s(0,t),V _{\sigma(0,t)}\rangle  \nonumber \\
&=& \frac{1}{r(t)}\langle \overline{\gamma}(t),V _{p}\rangle \equiv - \frac{a(t)}{r(t)}. \label{eq22}
\end{eqnarray}
Therefore, for all $t\in J$ for which $|N(t)|_h, |Y(t)|_h>0$,
\begin{eqnarray}
\left|h(\frac{N(t)}{|N(t)|_h}, \frac{Y(t)}{|Y(t)|_h})\right|  &=&  2 \frac{|\langle N(t), V_{\gamma(t)}\rangle |}{|N(t)|_h}  \frac{|\langle Y(t), V_{\gamma(t)}\rangle |}{|Y(t)|_h}\nonumber \\
&\stackrel{Eq.(\ref{eqclear})}{=}& \frac{2}{|N(t)|_h} |\langle N(t), V_{\gamma(t)}\rangle |  \frac{|h( Y(t), V_{\gamma(t)}) |}{|Y(t)|_h} \nonumber \\
&\stackrel{Eq.(\ref{eq22})}{\leq}& \frac{2|a(t)|}{r(t)|N(t)|_h} \nonumber \\
&\stackrel{Eqs.(\ref{eq21}),(\ref{eq22})}{=}&  \frac{2|a(t)|}{r(t) \sqrt{2\frac{a(t)^2}{r(t)^2} + 1}} \nonumber \\
&\equiv& \frac{2|a(t)|}{ \sqrt{2 a(t)^2 +r(t)^2}} \stackrel{Eq.(\ref{clearevenmore})}{\leq} \frac{2 \ell}{r(t)}, \label{eq23}
\end{eqnarray}
or, in other words, 
\begin{equation}\label{eq24}
|h(N(t),Y(t))|\leq \frac{2 \ell  |N(t)|_h|Y(t)|_h}{r(t)}. 
\end{equation} 
Let 
\[
\mathcal{J} := \{ t\in J \, : \, r(t) > \sqrt{8}\ell \}. 
\]
Since $r\equiv |\overline{\gamma}|$ is unbounded, this is not empty. Thus, for $t\in \mathcal{J}$, the inequality in (\ref{eq24}) yields
\begin{equation}\label{eq25}
|h(N(t),Y(t))|\leq \frac{1}{\sqrt{2}} |N(t)|_h|Y(t)|_h.
\end{equation} 
Finally, we compute, on  $ \mathcal{J}$,
\begin{eqnarray}
h(\gamma', \gamma') &=& (r')^2|N|^2_h + r^2|Y|^2_h + 2rr'h(N,Y) \nonumber \\
&\stackrel{Eq.(\ref{eq25})}{\geq}&   (r')^2|N|^2_h + r^2|Y|^2_h - 2r|r'||Y|_h|N|_h/\sqrt{2} \nonumber \\
&=& \frac{(r')^2|N|^2_h}{2} + \left( |r'||N|_h/\sqrt{2} - r|Y|_h\right)^2 \nonumber \\
& \stackrel{Eq.(\ref{eq21})}{\geq}& \frac{(r')^2}{2}. \label{eq26}
\end{eqnarray}
Now, note that since $r$ is still unbounded on $\mathcal{J}$, we can pick a sequence $(t_n)\subset \mathcal{J}$ for which $(r(t_n))$ is increasing and
\[
r(t_n) >\sqrt{2}n, \quad \forall n\in \mathbb{N}. 
\]
Therefore, integrating ``the square root'' of (\ref{eq26}), we obtain 
\[
L_h(\gamma) =\ell \geq n ,  \quad \forall n\in \mathbb{N},
\]
a contradiction. This concludes the proof. 
\qed  
\begin{cor}\label{stationary}
Let $(M,g)$ be a 
compact Lorentz manifold, and suppose there exists a timelike Killing vector field $V\in \mathfrak{X}(M)$ for which ${\rm Ric}(V,V)\leq 0$, where ${\rm Ric}$ denotes the Ricci tensor of $(M,g)$. Then, $\forall p \in M$ the exponential map $\exp_p:T_pM\rightarrow M$ is weakly proper. 
\end{cor}
\noindent {\em Proof.} From the hypotheses, $(M,g)$ must be geodesically complete (see \cite{RoSa}). Then, it is enough to use a standard Bochner-like trick to show that $V$ is parallel. See \cite[Thm. \ 3.2]{RoSan}\footnote{In \cite{RoSan} the authors actually require that $(M,g)$ be Ricci-flat, but their integral formula is easily seen to apply when one only has ${\rm Ric}(V,V)\leq 0$.}.
\qcd

\begin{notation}
Henceforth, we fix a connected affine manifold $M=( M,\nabla)$, and denote as before its pointwise exponential maps by $\exp_p:\mathcal{D}_p\subset T_pM\rightarrow M$, for every $p\in M$. Also, given each such $p\in M$, we denote the singularities of $\exp_p$ by
\[
\mathcal{S}_p:=\{ v\in \mathcal{D}_p \, : \, (d\exp _p)_v \mbox{ is singular}\},
\]
so that 
\[
{\rm Conj}(p):= \exp_p(\mathcal{S}_p)
\]
is the set of conjugate points to $p$ along some geodesic issuing from it. 
 Finally, we say that $(M,\nabla)$ is {\em weakly Wierdersehen} ({\em WW}) if
\begin{itemize}
\item[i)] ${\rm Conj}(p)$ is closed, and
\item[ii)] $M\setminus {\rm Conj}(p)$ is connected for all  $p \in M$.
\end{itemize}
\end{notation}

\bigskip
The round sphere $\mathbb{S}^n$ is an obvious example of a WW manifold. Of course, it is also geodesically connected by Hopf-Rinow, and its exponential maps at each point are weakly proper at each point by Prop. \ref{prop3}. However, the following example shows that in more general semi-Riemannian contexts  some condition along the lines of weak properness is likely needed if geodesic connectedness is to be expected, even for WW manifolds. 
\begin{exe}\label{deSitter}
Consider, for a fixed $n \in \mathbb{N}$, the $(n+1)$-dimensional {\em Minkowski spacetime}, i.e., $\mathbb{R}^{n+1}$ with the flat Lorentzian metric given in standard Cartesian coordinates $(x_1, \ldots, x_{n+1})$ by 
\[
\eta ^{n+1}= \sum ^{n}_{i=1} dx_i^2 - dx_{n+1}^2.
\]
We write $\mathbb{R}^{n+1}_1= (\mathbb{R}^{n+1}, \eta^{n+1})$. The {\em de Sitter space} $\mathbb{S}^n_1$ is then the one-sheeted hyperboloid  
\[
 \sum ^{n}_{i=1} x_i^2 - x_{n+1}^2 =1
 \]
endowed with the Lorentzian metric induced by $\eta^{n+1}$. This is a geodesically complete $n$-dimensional Lorentzian manifold of constant sectional curvature 1, and hence can be regarded as the Lorentzian analogue of the round sphere. The following facts about $\mathbb{S}^n_1$ are well-known (cf., e.g., Prop. 4.28 in \cite{oneill}). 
\begin{itemize}
\item[1)] There are no conjugate points along either lightlike or timelike geodesics, but there are conjugate points along spacelike geodesics.
\item[2)] Up to affine reparametrization, any spacelike geodesic $\alpha:\mathbb{R}\rightarrow \mathbb{S}^n_1$ is given as follows. Let $\{e_1,e_2\}$ be a pair of unit spacelike vectors in $\mathbb{R}^{n+1}_1$, orthogonal to one another. Then 
\begin{equation}\label{spacegeod}
\alpha(t) = \cos t \cdot e_1+ \sin t\cdot e_2, \quad \forall t \in \mathbb{R}. 
\end{equation}
\end{itemize}
Let $p = \alpha(0)=e_1$. Clearly, $T_p\mathbb{S}^n_1$ is essentially the orthogonal complement $e_1^{\perp}$ of $e_1$ in $\mathbb{R}^{n+1}_1$, and the only conjugate points to $p$ are $p$ itself and $-p$($\equiv \alpha(k\pi)$, for any odd $k\in \mathbb{Z}$). In particular  ${\rm Conj}(p) =\{p,-p\}$, and thus $\mathbb{S}^n_1$ is a WW manifold. The form (\ref{spacegeod}) then shows that one can identify $\mathcal{S}_p$ with the set of all spacelike vectors in $T_p\mathbb{S}_1^n\simeq \mathbb{R}^n_1$ of norm $k \pi$, for $k \in \mathbb{N}$: if we use Cartesian coordinates $(z_1, \ldots, z_n)$ therein, any connected component of this set is diffeomorphic to the hyperboloid 
\[
 \sum ^{n}_{i=1} z_i^2 - z_{n+1}^2 =1.
 \]
Since the latter set is non-compact, $\exp_p$ cannot be weakly proper. Surely enough,  $\mathbb{S}^n_1$ is not geodesically connected. 
\end{exe}
\section{Proof of Theorem \ref{mainthm}}\label{s2}

Theorem \ref{mainthm} follows immediately from the following more general result, by taking $F=\exp _p$, $N=\mathcal{D}_p\subset T_pM$ its maximal domain, and $O=0_p \in \mathcal{D}_p$ the zero vector.
\begin{thm}\label{mmain}
Let $F:N\rightarrow M$ be a smooth map between smooth manifolds of same dimension, which is regular and weakly proper at $O\in N$. Denote by $\mathcal{S}$ the set of singular points of $F$. If $F(\mathcal{S})$ is closed and $M\setminus F(\mathcal{S})$ connected, then $F$ is surjective.
\end{thm}
In order to prove that $F:N\rightarrow M$ is surjective, denote $p=F(O)$ and fix any other $q \in M$. There is no loss of generality in taking $q\in M\setminus F(\mathcal{S})$. 
From the hypotheses, $\mathcal{R}:= M\setminus F(\mathcal{S})$ is an open connected subset of $M$ containing $q$. Take some point $r\in \mathcal{R}\cap \mathcal{V}$, where $\mathcal{V}=F(\mathcal{U})\ni p$ is an open set such that $U$ is connected and $F|_U:U\subset N \rightarrow \mathcal{V}$ is a diffeomorphism\footnote{When $F=\exp _p$ we chose $\mathcal{V}$ to be a normal neighborhood of $p$.} (note that $F(\mathcal{S})$ has zero measure - and hence empty interior - by Sard's theorem; hence $\mathcal{V}\not\subset F(\mathcal{S})$, i.e., $\mathcal{R}\cap \mathcal{V}$ is non-empty).

Since $\mathcal{R}$ is an open submanifold of $M$, we can fix a background complete Riemannian metric $h$ on $\mathcal{R}$. Recall also that $\mathcal{R}$ is connected and $r,q\in \mathcal{R}$.  
Pick an $h$-geodesic $\sigma:[0,l]\rightarrow \mathcal{R}$ connecting $r$ to $q$, 
parametrized by $h$-arc length. 
Let $\mathcal{C}:=F^{-1}(\mathcal{R})\subset N$. It is clear that $\mathcal{C}\subset N\setminus \mathcal{S}$ is open, and thus  $\phi:=F\mid_{\mathcal{C}}:\mathcal{C}\rightarrow \mathcal{R}$ is a local diffeomorphism. Thus, $\overline{h}:=(\phi)^*h$ is a Riemannian metric on $\mathcal{C}$ and the map 
\[
\phi: (\mathcal{C},\overline{h})\rightarrow (\mathcal{R},h)
\]
is a local {\em isometry}.

Since $r\in \mathcal{R}\cap \mathcal{V}$, in particular, $r$ is in the image of $F$ and there exists some $u\in \mathcal{C}(=F^{-1}(\mathcal{R}))$ such that $\phi(u)=r$. 
Let $\overline{\sigma}:[0,\overline{\ell})\rightarrow \mathcal{C}$ be the right-inextendible $\overline{h}$-geodesic such that $\overline{\sigma}(0)=u$ and $(d \phi)_{u}(\overline{\sigma}'(0))\equiv \sigma'(0)$ ($\overline{\sigma}'(0)$ exists because $\phi$ is a local diffeomorphism). Since $\phi: (\mathcal{C},\overline{h})\rightarrow (\mathcal{R},h)$
is a local isometry, $\phi\circ \overline{\sigma}$ is an $h$-geodesic, and it follows from uniqueness that  
\[
\phi\circ \overline{\sigma}|_{[0,\overline{\ell})\cap [0,\ell]} = \sigma |_{[0,\overline{\ell})\cap [0,\ell]}.
\]
Therefore, if $\ell < \overline{\ell}$, we would have
\[
F(\overline{\sigma}(\ell))=\phi\circ\overline{\sigma}(l)=\sigma(\ell)=q,
\]
i.e., $q\in {\rm Im}(F)$. Thus, suppose by contradiction that  $\overline{\ell}\leq \ell$. Then $F\circ \overline{\sigma} \equiv \sigma |_{[0,\overline{\ell})}$ is right-extendible, and by weak properness, there exists a compact set $K\subset N$ containing the image of $\overline{\sigma}$. 
Now, let $(t_k),(s_k)\subset [0,\overline{\ell})$ be any two sequences, both converging to $\overline{\ell}$. On the one hand, the fact that $\overline{\sigma}(t_k),\overline{\sigma}(s_k)\subset K$ implies that up to passing to subsequence we can take
\[
\overline{\sigma}(t_k)\rightarrow v \mbox{  and  } \overline{\sigma}(s_k)\rightarrow w \quad v,w \in N.
\]
On the other hand, $F(\overline{\sigma}(t_k)),F(\overline{\sigma}(s_k))\rightarrow \sigma(\overline{\ell}) \in \mathcal{R}$. We conclude that $v,w\in \mathcal{C}$ and 
\[
F(v)=F(w).
\]
However, $\overline{\sigma}$ is a unit $\overline{h}$-geodesic, so
\[
d_{\overline{h}}(\overline{\sigma}(t_k),\overline{\sigma}(s_k))\leq |t_k-s_k|\rightarrow 0,
\]
whence we conclude that $v=w$. But we have now proven that $v$ is a right-endpoint of $\overline{\sigma}$, in contradiction. Thus, $\ell <\overline{\ell}$ and the proof is complete.

\section*{Acknowledgments}
The authors are partially supported by the project MTM2016-78807-C2-2-P (Spanish MINECO with FEDER funds). 

\end{document}